\documentclass[10pt]{article}
\usepackage{mathrsfs}
\usepackage{amsmath}
\usepackage{amssymb}
\usepackage{graphicx}
\def \P{\mathrm P}
\newtheorem{theorem}{\scshape \mdseries  Theorem}[section]
\newtheorem{lemma}[theorem]{\scshape \mdseries  Lemma}
\newtheorem{coro}[theorem]{\scshape \mdseries  Corollary}

\newtheorem{rem}[theorem]{\scshape \mdseries  Remark}

\begin{document}

\title{\sf The Connectivity and the Harary Index of a Graph\thanks{
      Supported by National Natural Science Foundation of China (11071002),
Program for New Century Excellent Talents in University,
Key Project of Chinese Ministry of Education (210091),
Specialized Research Fund for the Doctoral Program of Higher Education (20103401110002),
Science and Technological Fund of Anhui Province for Outstanding Youth  (10040606Y33),
Scientific Research Fund for Fostering Distinguished Young Scholars of Anhui University(KJJQ1001),
 Academic Innovation Team of Anhui University Project (KJTD001B). 
   }}
\author{Xiao-Xin Li$^{1,2}$, Yi-Zheng Fan$^{2,}$\thanks{Corresponding author. E-mail addresses:
fanyz@ahu.edu.cn (Y.-Z. Fan), lxx@czu.edu.com (X.-X. Li)}
\\
  {\small  \it $^1$Department of Mathematics, Chizhou University, Chizhou 247000, P. R. China} \\
  {\small \it $^2$School of Mathematical Sciences, Anhui University, Hefei 230601, P. R. China}
 }
\date{}
\maketitle

\noindent {\bf Abstract}\ \ The Harary index of a graph is defined as the sum of reciprocals of distances between all pairs of vertices of the graph.
In this paper we provide an upper bound of the Harary index in terms of the vertex or edge connectivity of a graph.
We characterize the unique graph with maximum Harary index among all graphs with given number of cut vertices or vertex connectivity or edge connectivity.
In addition we also characterize the extremal graphs with the second maximum Harary index among the graphs with given vertex connectivity.

\noindent {\bf MR Subject Classifications:} 05C90


\section{Introduction}
Let $G$ be a simple graph with vertex set $V(G)$ and edge set $E(G)$.
The {\it distance} between two vertices $u,v $ of $G$, denoted by $d_G(u,v)$, is defined as the minimum length of the paths between $u$ and $v$ in $G$.
The Harary index of a graph $G$, denoted by $H(G)$, has been introduced
independently by Plav\v{s}i\'c et al. \cite{pla} and by Ivanciuc et al. \cite{iva} in 1993
for the characterization of molecular graphs. It has been named in honor of Professor Frank Harary on the occasion of his 70th
birthday.
 The {\it Harary index} $H(G)$ is defined as the sum of reciprocals of distances between all pairs of vertices of the graph $G$, i.e. 
 $$H(G)=\sum_{u,v \in V(G)} \frac{1}{d_G(u,v)}.$$ Mathematical properties and applications of the Harary index are reported in \cite{diu,est,he,lu,zhou}.
Note that in any disconnected graph $G$, the distance is infinite between any two vertices from two distinct components.
Therefore its reciprocal can be viewed as $0$. Thus, we can define validly the Harary index of disconnected graph $G$ as follows:
$$H(G)=\sum_{i=1}^k H(G_i),$$ where $G_1, G_2, \ldots, G_k$ are the components of $G$.

 Another distance-based topological index of a graph $G$ is the Wiener index, denoted by $W(G)$. As an oldest topological index, the {\it Wiener index} of a graph $G$, first introduced by Wiener \cite{wie} in 1947, was defined as  $$W(G)=\sum_{u,v \in V(G)} d_G(u,v).$$ The motivation for introduction of the Harary index was pragmatic - the aim was
to design a distance index differing from the Wiener index in that the contributions
to it from the distant atoms in a molecule should be much smaller than from near
atoms, since in many instances the distant atoms influence each other much less than
near atoms.

Let $\gamma (G,k)$ be the number of vertex pairs of the graph $G$ that are at distance $k$. Then
$$H(G)=\sum_{k\geq 1}\frac{1}{k} \gamma (G,k). \eqno (1.1)$$
It will be convenient to determine the exact value by Eq. (1.1) for some graphs with simple structure (e.g. the graphs with small diameter),
but in general it is very difficult to give the exact value of $\gamma (G,k)$.
So it is very useful to provide upper or lower bounds for the Harary index;  see e.g. \cite{das,he,zhou}.
In addition, the extremal Harary index of a given class of graphs has also been studied extensively; see e.g. \cite{feng,ili,xu1,xu2,xu3,xu4,yu}.

In this paper we provide an upper bound of the Harary index in terms of the vertex or edge connectivity of a graph.
We characterize the unique graph with maximum Harary index among all graphs with given number of cut vertices or vertex connectivity or edge connectivity.
In addition we also characterize the extremal graphs with the second maximum Harary index among the graphs with given vertex connectivity.

\section{Main results}
In Section 2.1 we determine the unique graph with maximum Harary index among all graphs with given number of cut vertices.
We find the optimal graph is surely connected with vertex or edge connectivity $1$.
In Section 2.2 we consider a general problem, that is, determining the graph(s) with maximum Harary index among all graphs with fixed vertex or edge connectivity.
By these results we provide an upper bound of Harary index of a graph in terms of the vertex or edge connectivity.

We introduce some notions used in this paper.
Let $G$ be a graph.
For a vertex $v \in V(G)$, denote by $N_G(v)$ the neighborhood of $v$ in $G$ and by $d_G(v)=|N_G(v)|$ the degree of $v$ in $G$.
A vertex of $G$ is called {\it pendent}  if it has degree $1$, and the edge incident with a pendent vertex is a {\it pendent edge}.
A {\it pendent path} at $v$ in a graph $G$ is a path in which no vertex other than $v$ is incident with any edge of $G$ outside the path, 
where the degree of $v$ is at least three.
 %In particular, if a vertex $v$ is neither a pendent vertex nor a cut vertex of $G$, then it is considered as a pendant path of length $0$ at itself.
A {\it cut vertex} of a graph is a vertex whose removal increases the number of components of the graph.
A {\it block} of a connected graph is defined to be a maximum connected subgraph without cut vertices.
The {\it vertex connectivity} (respectively, {\it edge connectivity}) of a graph, 
is the minimum number of vertices (respectively, minimum number of edges) whose deletion yields the resulting graph disconnected or a singleton.
%For a complete graph $K_n$, we define by convention that $\kappa(K_n)=n-1$.

For a subset $W \subset V(G)$, let $G-W$ be the subgraph of $G$ obtained by deleting the vertices of $W$ together with the edges incident with them.
Similarly, for a subset $E_1 \subset E(G)$, denote by $G-E_1$ the subgraph of $G$ obtained by deleting the edges of $E_1$.
For an edge set $E_2 \nsubseteq E(G)$, if two endpoints of any edge in $E_2$ belong to $V(G)$, then we denote by $G+E_2$ the graph obtained from $G$ by adding the edges of $E_2$.
Denote by $P_n=\P v_1v_2\cdots v_n$ a path on vertices $v_1,v_2,\ldots,v_n$ with edges $v_iv_{i+1}$ for $i=1,2,\ldots,n-1$;
and denote by $K_n$ the complete graph on $n$ vertices.

%According to Lemma 2.1, for a graph $G$ on $n$ vertices, we have $0 \leq H(G) \leq \frac{n(n-1)}{2}$, with left equality holds if and only if $G$ is a graph with no edge and right equality holds if and only if $G=K_n$.
%
%For a graph $G$ with $v_i \in V(G)$, Xu and Ch. Das \cite{xu2} defined a function $Q_G(v_i)=\sum_{v_j \in V(G)}\frac{d_G(v_i,v_j)}{d_G(v_i,v_j)+1}.$ And with the help of this function, they give two lemmas on the effect on the Harary index of a graph by "grafting" an edge. Note that Lemma  is also presented in \cite{he}, which is proved by another method.
%
%\begin{lemma}\cite{he,xu2}
%Let $G$ be a connected graph on $n\geq 2$ vertices and $u$ be a vertex of $G$. Let $G_{k,l}$ be the graph obtained from $G$ by attaching two new paths $P:uv_1v_2\ldots v_k$ and $Q:uu_1u_2\ldots u_l$ of length $k$ and $l$ at $u$, respectively, where $u_1,\ldots,u_l$ and $v_1,\ldots,v_k$ are distinct new vertices. If $k\geq l\geq 1$, then $$H(G_{k,l})>H(G_{k+1,l-1}).$$
% \end{lemma}
%
%\begin{lemma}\cite{xu2}
%Let $G$ be a connected graph on $n > 2$ vertices and $v_i,v_j$ be its distinct vertices with $Q_G(v_i)=Q_G(v_j)$. Let $G_{s,t}^*$ be the graph obtained from $G$ by attaching at $v_i$ one path $P:v_iu_1u_2\ldots u_s$ and at $v_j$ the other path $Q:v_ju_1'u_2'\ldots u_t'$ of length $s$ and $t$, respectively, where $u_1,\ldots,u_s$ and $u_1',\ldots,u_t'$ are distinct new vertices. If $s\geq t\geq 1$, then $$H(G_{s,t}^*)>H(G_{s+1,t-1}^*).$$
% \end{lemma}

\subsection{Maximum Harary index with given number of cut vertices}
%For two vertex-disjoint connected graphs $G_1$ and $G_2$ with $uv\in E(G_1)$ and $x_0\in V(G_2)$, and for a positive integer $s$, let $G_1uG_2(s)$ be the graph obtained from $G_1$ and $G_2$ by identifying $u$ and $x_0$, and adding a path $P_s=z_1z_2\ldots z_s$ with $z_1=v$.

\begin{lemma}\cite{xu2} \label{addedge}
Let $G$ be a graph with $u,v \in V(G)$. If $uv \notin E(G)$, then $H(G) < H(G+uv)$. If $uv \in E(G)$, then $H(G) > H(G-uv)$.
 \end{lemma}

\begin{lemma} \label{per1} Let $G_1,G_2,P_s$ be pairwise vertex-disjoint connected graphs,
where $G_1$ contains an edge $uv$ such that $N_{G_1}(u)\setminus \{v\}=N_{G_1}(v)\setminus \{u\}=\{w_1, w_2, \ldots, w_k\}\;(k \ge 1)$,
$G_2$ contains a shortest path $\P x_1\cdots x_t\;(t \ge s+2) $ from $x_1$ to $x_{t}$, and $P_s=\P z_1z_2 \ldots z_s$.
Let $G$ be obtained from $G_1$ by identifying $u$ with $x_1$ of $G_2$ and identifying $v$ with $z_1$ of $P_s$,
and let $G'=G-\{vw_1,vw_2,\ldots,vw_k\}+\{x_2w_1,x_2w_2,\ldots,x_2w_k\}.$
Then $$H(G)<H(G'),$$
where the graphs $G$ and $G'$ are shown in Fig. 2.1.
\end{lemma}

{\it Proof:}
Let $P$ be the path of $G$ obtained by connecting the paths $\P x_1\cdots x_t$, $\P uv$ and $\P z_1z_2 \ldots z_s$,
where $u=x_1$ and $v=z_1$.
Partition the vertex set of $G$ as
$$V(G)=(V(G_1)\backslash\{u,v\}) \cup (V(G_2)\backslash\{x_i: i=1,2,t\}) \cup V(P)=:S_1 \cup S_2 \cup S_3.$$
From $G$ to $G'$, the distance between any two vertices in each $S_i$ is unchanged for $i=1,2,3$;
the distance from any vertex of $S_1$ to any of $S_2$ is not increased;
the distance from any vertex of $S_1$ to any of $z_i \;(i=1,\ldots,s)$ of $S_3$ is increased by $1$, and to any of $x_i\;(i=2,\ldots,t)$ are decreased by $1$,
and to the vertex $u$ is unchanged;
the distance from any vertex of $S_2$  and any of $S_3$ is unchanged.

For any vertex $y \in S_1$, assuming that $d_{G'}(y,z_1)=a \;( \ge 2)$, then
$d_{G}(y,z_1)=a-1, d_{G'}(y,x_2)=a-1, d_{G}(y,x_2)=a$.
Thus
\begin{align*}
\Delta(y) & := \sum \limits_{i=1}^s\frac{1}{d_{G''}(y,z_i)}+\sum \limits_{i=2}^{t}\frac{1}{d_{G''}(y,x_i)}
-\sum \limits_{i=1}^s\frac{1}{d_{G'}(y,z_i)}-\sum \limits_{i=2}^{t}\frac{1}{d_{G'}(y,x_i)}   \\
& =  \sum \limits_{i=0}^{s-1}\frac{1}{a+i}+\sum \limits_{i=0}^{t-2}\frac{1}{a-1+i}-\sum \limits_{i=0}^{s-1}\frac{1}{a-1+i}-\sum \limits_{i=0}^{t-2}\frac{1}{a+i}\\
 & = \sum \limits_{i=s}^{t-2}\left(\frac{1}{a-1+i}-\frac{1}{a+i}\right)\\
 & >  0.
 \end{align*}
So $$H(G')-H(G)\geq \sum \limits_{y \in S_1} \Delta(y)>0.$$
The result follows.\hfill$\blacksquare$

\begin{center}
\vspace{3mm}
\includegraphics[scale=.6]{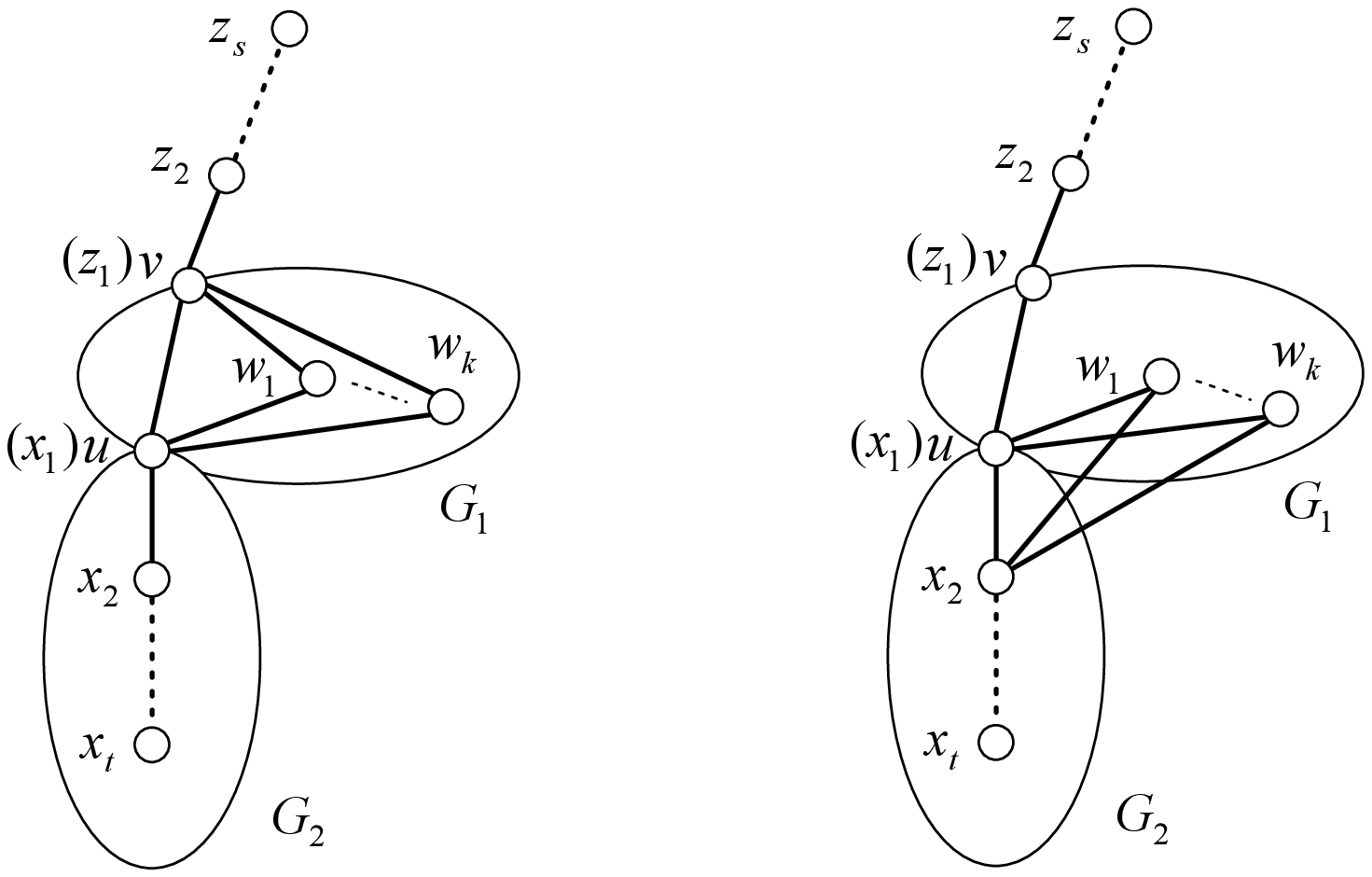}
\vspace{3mm}

{\small Fig. 2.1. The graphs $G$ (left) and $G'$ (right) in Lemma \ref{per1}.}

\end{center}

\begin{rem} \label{rem1}
The graphs $G$ and $G'$ in Lemma \ref{per1} possess the same number of cut vertices.
In addition, if taking $s=1$ (i.e. the path attaching at $v$ is trivial), the edge $uv$ of $G$ will become a pendant edge of $G'$.
\end{rem}

If taking $G_2=\P x_1\cdots x_t$ in Lemma \ref{per1}, we will have the result below, which has been shown in \cite{xu2} and \cite{he} under a more general condition. 

\begin{coro} \label{coro-per1}
Let $G $ be a connected graph containing an edge $uv$ such that $N_{G}(u)\setminus \{v\}=N_{G}(v)\setminus \{u\} \ne \emptyset$.
Let $G_{t,s}$ be obtained from $G$ by attaching a path $P_t$ at $u$ and a path $P_s$ at $v$.
If $t \ge s+2 \ge 3$, then $H(G_{t,s}) < H(G_{t-1,s+1})$.
\end{coro}

%For two vertex-disjoint complete graphs $K_p, K_q (p,q>2)$, denote by $K_puK_q$ the graph obtained from $K_p$ and $K_q$ by identifying one vertex of $K_p$ and one vertex of $K_q$, where the unique common vertex of $K_p$ and $K_q$ is denoted by $u$. Except for $u$, attaching a path $P_s (s>1)$ on each vertex of $K_q$ and one vertex of $K_p$, attaching other components on the other vertices of $K_p$, all these graphs consist of a set, denoted by $\mathcal{G}(p,q,s)$.

%\begin{lemma}
%Let $K_puK_q$ be the union of two complete graphs $K_p,K_q$ sharing exactly one common vertex $u$.
%Let $G$ be obtained from $K_puK_q$ by attaching a path $P_s \; (s >1)$ at each vertex of $V(K_q)\setminus\{u\}$ and some vertex $w_1 \in V(K_p)\setminus\{u\}$,
%and attaching some connected graph at each of $V(K_p)\setminus\{u,w_1\}$;
%and let $G'$ be obtained from $G$ by deleting the edges of $K_q$ incident to $v_1$ except $v_1u$, and adding all possible edges between each of $V(K_q)\setminus\{v_1\}$ and each of $V(K_p)$; see Fig. 2.2 for the graph $G$ and $G'$.
%Then $H(G)<H(G')$.
%\end{lemma}

\begin{lemma} \label{per2}
Let $K_puK_q$ be the union of two complete graphs $K_p,K_q$ sharing exactly one common vertex $u$, where $p \ge 3, q \ge 3$.
Let $G$ be obtained from $K_puK_q$ by attaching a path $P_t$ at some vertex $w_1 \in V(K_p)\backslash\{u\}$ and a path $P_s$ at some vertex $v_1 \in V(K_q)\backslash\{u\}$,
and possibly attaching some connected graphs at other vertices of $V(K_puK_q)\backslash\{u,v_1,w_1\}$, where $t \ge s \ge 1$; and
let $G'$ be obtained from $G$ by deleting the edges of $K_q$ incident to $v_1$ except $v_1u$ and adding all possible edges between each of $V(K_q)\backslash\{v_1\}$ and each of $V(K_p)$; see Fig. 2.2 for the graph $G$ and $G'$.
Then $H(G)<H(G')$.
\end{lemma}

{\it Proof:}
Let $G^*$ be the component of $G-\{u,v_1\}$ which contains the vertices of $K_q$,
and $G^{**}$ be the component of $G-\{u,w_1\}$ which contains the vertices of $K_p$.
Let $P_s=\P v_1v_2\cdots v_s$, $P_t=\P w_1w_2\cdots w_t$.
Partition the vertex of $G$ as
$$V(G^*) \cup V(P_s) \cup \{u\} \cup V(P_t) \cup V(G^{**})=:S_1 \cup S_2 \cup S_3 \cup S_4 \cup S_5.$$
Observe the transformation from $G$ to $G'$, the distances between any vertex of $S_1$ and any of $S_2$ are increased by $1$,
the distances between any vertex of $S_1$ and $S_4$ are decreased by $1$,
the distance between any vertex of $S_1$ and any of $S_5$ is decreased by $1$,
and the distance between any other two vertices is not changed.

For any vertex $y \in S_1$, assuming that $d_{G'}(y,v_1)=a (\ge 2)$, then
$d_{G'}(y,w_1)=a-1, d_{G}(y,v_1)=a-1, d_{G}(y,w_1)=a.$
Thus
\begin{align*}
\Delta(y) &:= \sum \limits_{i=1}^s\frac{1}{d_{G'}(y,v_i)}+\sum \limits_{i=1}^t\frac{1}{d_{G'}(y,w_i)}
-\sum \limits_{i=1}^s\frac{1}{d_{G}(y,v_i)}-\sum \limits_{i=1}^t\frac{1}{d_{G}(y,w_i)}   \\
& =  \sum \limits_{i=0}^{s-1}\frac{1}{a+i}+\sum \limits_{i=0}^{t-1}\frac{1}{a-1+i}-\sum \limits_{i=0}^{s-1}\frac{1}{a-1+i}-\sum \limits_{i=0}^{t-1}\frac{1}{a+i}\\
 & =  \sum \limits_{i=s}^{t-1}\left(\frac{1}{a-1+i}-\frac{1}{a+i}\right) \ge 0.
 \end{align*}
So, $$H(G')-H(G)= \sum \limits_{y \in S_1}\Delta(y) + \sum \limits_{(y,z) \in S_1 \times S_5}\left(\frac{1}{d_{G'}(y,z)}-\frac{1}{d_G(y,z)}\right)
>\sum \limits_{y \in S_1}\Delta(y) \ge 0.$$
The result follows.\hfill$\blacksquare$

\begin{center}
\vspace{3mm}
\includegraphics[scale=.6]{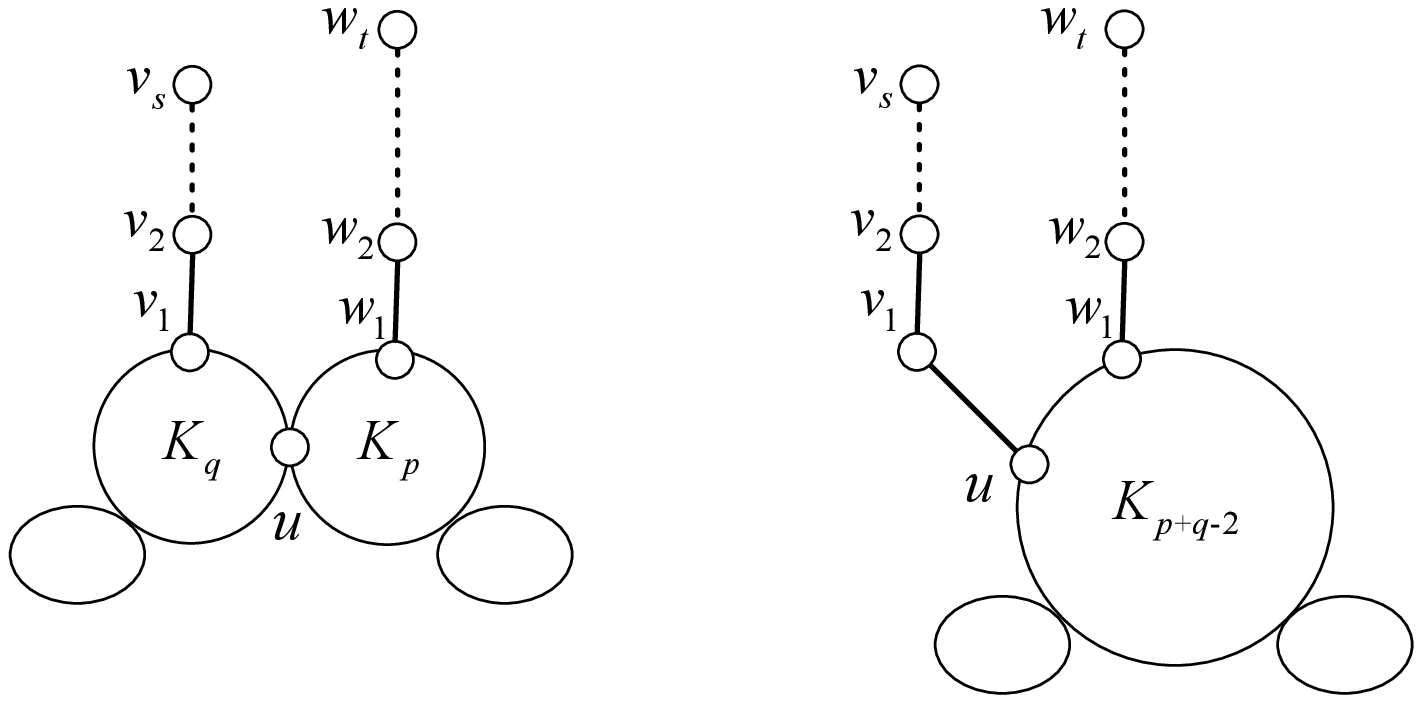}
\vspace{3mm}

{\small Fig. 2.2.  The graphs $G$ (left) and $G'$ (right) in Lemma \ref{per2}}
\end{center}

\begin{rem} \label{rem2}
The graphs $G$ and $G'$ in Lemma \ref{per2} possess the same number of cut vertices.
%Furthermore, if the graphs $G^*,G^{**}$ both have at least one vertices, then the result also holds without the limitation of $p \ge 3,q \ge 3$.
In addition, if taking $s=1$, the edge $uv_1$ of $G$ becomes a pendant edge of $G'$.
\end{rem}

\begin{theorem} \label{cutvertex}
Among all graphs with $n$ vertices and $k$ cut vertices, where $0\leq k \leq n-2$,
the maximal Harary index is attained uniquely at the graph $\mathbf{G}_{n,k}$,
where $\mathbf{G}_{n,k}$ is obtained from $K_{n-k}$ by attaching $n-k$ paths of almost equal lengths to its vertices respectively.
\end{theorem}

{\it Proof:} Let $G$ be a graph with the maximal Harary index among all the graphs with $n$ vertices and $k$ cut vertices.
If $k=0$, then by Lemma \ref{addedge}, $G = K_n = \mathbf{G}_{n,0}$. Suppose in the following that $1\leq k \leq n-2$.
The result will hold by the following claims.

{\it Claim 1: $G$ is connected.} Assume that $G$ is disconnected.
Let $z$ be a cut vertex of $G$. Then $z$ is also a cut vertex of some component, say $G_1$, of $G$.
Let $G_2$ be a component of $G$ different from $G_1$. If there is a cut vertex, say $z'$, in $G_2$, then $G+zz'$ possesses the same number of cut vertices as $G$,
and by Lemma \ref{addedge}, $H(G)<H(G+{zz'})$, a contradiction.
If there are no cut vertices in $G_2$, adding edges between $z$ and all vertices of $G_2$, we will arrive at a new graph $G'$, which
possesses the same number of cut vertices as $G$.
However, by Lemma \ref{addedge}, $H(G)<H(G')$, a contradiction again. So $G$ is connected.

By Lemma \ref{addedge}, each block of $G$ is complete, and each cut vertex of $G$ is contained in exactly two blocks.
If every block of $G$ is trivial (containing exactly two vertices), i.e., every block is a single edge, then $G$ is a tree with maximum degree two, i.e., $G = P_n= \mathbf{G}_{n,n-2}$.
Suppose in the following that $G$ contains nontrivial blocks (on at least three vertices).

{\it Claim 2: If $G \ne \mathbf{G}_{n,1}$, then each pendent block (i.e., the bolck containing only one cut vertex of $G$) is an edge.}
Assume to the contrary, $B_1$ is a nontrivial pendent block of $G$.
Let $u$ be a vertex of $B_1$ different from the unique cut vertex, say $w$, contained in $B_1$.
Let $B_2$ be the block adjacent to $B_1$. Deleting all edges in $B_1$ incident to $u$ except $uw$,
and adding all edges between the vertices of $V(B_1)\backslash \{u\}$ and the vertices of $V(B_2)$, we obtain a new graph $G'$ with the same number of cut vertices as $G$.
By Remark \ref{rem1} and Remark \ref{rem2}, and the fact that $G \ne \mathbf{G}_{n,1}$, we have $H(G)<H(G')$, a contradiction.

Choose a pendent path, say $P_s$ attached at $v$ of some nontrivial block $B$, whose length is minimum among all pendant paths of $G$.
We stress that $P_s$ may be trivial, i.e. $s=1$ or $P_s$ contains only the vertex $v$.

{\it Claim 3: The component attached at any vertex of $B$ is a path (possibly being trivial)}.
For $x\in V(B)$, let $H^{(x)}$ be the component of $G-E(B)$ containing $x$.
Obviously, $H^{(v)}= P_s$. Suppose $u$ is an arbitrary  vertex of $B$ and $u \neq v$.
Obviously, $N_B(v)\backslash \{u\}=N_B(u)\backslash \{v\}$.
Let $G_1$ be the component of $G-E(H^{(u)}\cup E(P_s)$ containing $u$, which surely contains the block $B$.
%Then $G\cong G^*uH^{(u)}(s)$.

Assume that $H^{(u)}$ is not a (possibly trivial) path. Then $H^{(u)}$ contains a nontrivial block.
By the proof of Claim 2, $H^{(u)}$ must contain a nontrivial pendant path $P_t$ attached at some nontrivial block $B'$ of $H^{(u)}$, where $t \ge s$.
So $H^{(u)}$ contains a shortest path $P_r$ from $u$ to the pendant vertex of $P_t$, where $r \ge 3$ and $r \ge t+1 \ge s+1$.
If $r \ge s+2$ or $s=1$, by Lemma \ref{per1} and Remark \ref{rem1}, we may get another graph with $n$ vertex and $k$ cut vertices, which has a larger Harary index, a contradiction.
So, it suffices to consider the case: $s>1$, $B'$ shares with $G_1$ (also the block $B$) the common vertex $u$, and $H^{(u)}$ is obtained from $B'$ by attaching $P_s$ at each of its vertices except $u$.
Now applying Lemma \ref{per2}, we may get a new graph with $n$ vertices and $k$ cut vertices, which has a larger Harary index, a contradiction.
Hence $H^{(u)}$ is a pendent path attached at $u$ which contains at least $s$ vertices.

{\it Claim 4: All paths attached at the vertices of $B$ have almost equal lengths.}
This can be shown by Corollary \ref{coro-per1}.
\hfill$\blacksquare$

\begin{theorem}\label{cut1}
$$ H(K_n)=H(\mathbf{G}_{n,0}) > H(\mathbf{G}_{n,1})>H(\mathbf{G}_{n,2})> \cdots > H(\mathbf{G}_{n,n-2})=H(P_n).$$
Furthermore, if a graph $G$ of order $n \ge 3$ contains cut vertices or cut edges, then
$$H(G) \le H(\mathbf{G}_{n,1}),$$
with equality if and only if $G=\mathbf{G}_{n,1}$.
\end{theorem}

{\it Proof:}
Let $G=\mathbf{G}_{n,k}$, where $k \ge 1$.
Let $P_s$ be a pendant path of $G$ attached at $u$ with maximum length.
Surely $s \ge 2$.
Let $v$ be the vertex on the path $P_s$ adjacent to $u$.
Adding all possible edges between $v$ and the vertices of $K_{n-k}$ (the subgraph of $G$),
we will arrive at a graph $G' = \mathbf{G}_{n,k-1}$, which holds that $H(\mathbf{G}_{n,k-1})>H(G)=H(\mathbf{G}_{n,k})$ by Lemma \ref{addedge}.
The first assertion follows.
The remaining parts of this theorem can be obtained from above discussion and Theorem \ref{cutvertex}. \hfill$\blacksquare$

\subsection{Maximum Harary indices with given connectivity and edge connectivity}
In Section 2.1 we have determined the unique graph with maximum Harary index among all graphs with vertex or edge connectivity $1$; see Theorem \ref{cut1}.
Now we consider a general problem, i.e. characterizing the graph(s) with maximum Harary index among all graphs with fixed vertex or edge connectivity $k$.

We first give some notations.
For two vertex-disjoint graphs $G$ and $H$, let $G \cup H$ denote the union of $G$ and $H$, and $G \vee H$ denote the graph obtained from $G \cup H$
by adding all possible edges between the vertices of $G$ and the vertices of $H$.
Denote $G_{n_1,n_2,n_3}:=(K_{n_1}\cup K_{n_2})\vee K_{n_3},$ where $n_1\geq n_2 \geq 1$ and $n_3 \geq 1$.
Denote by $\mathcal{G}_n^r$ (respectively, $\overline{\mathcal{G}}_n^r$) the set of all connected graphs of order $n$ with vertex connectivity $r$ (respectively,  edge connectivity $r$).
Clearly, $1\leq r \leq n-1$, and $\mathcal{G}_n^{n-1}=\overline{\mathcal{G}}_n^{n-1}=\{K_n\}$.
So it is enough to consider the case of $1 \le r \le n-2$.
 Let $K(n-1,r)$ be a graph obtained from $K_{n-1}$ by adding a vertex together with edges joining this vertex to $r$ vertices of $K_{n-1}$,
 where  $1 \le r \le n-2$.
 Surely $K(n-1,r) \in \mathcal{G}_n^r$ and $K(n-1,r)\in \overline{\mathcal{G}}_n^r$.

\begin{lemma} \label{con}
If $n_1 \geq n_2 \geq 2$ and $n_3 \geq 1$, then $H(G_{n_1,n_2,n_3})< H(G_{n_1+1,n_2-1,n_3})$.
\end{lemma}

{\it Proof:}
Observe that the graph $G_{n_1,n_2,n_3}$ can be considered as one obtained from $G_{n_1,n_2-1,n_3}$ by adding a vertex, say $u$, and connecting $u$ with all vertices of $K_{n_2-1}\cup K_{n_3}$, and $G_{n_1+1,n_2-1,n_3}$ can be considered as one also obtained from $G_{n_1,n_2-1,n_3}$ by adding a vertex, also say $u$ for simplicity, and connecting $u$ with all vertices of $K_{n_1}\cup K_{n_3}$.
So, from $G_{n_1,n_2,n_3}$ to $G_{n_1+1,n_2-1,n_3}$, the distance between $u$ and any vertex of $K_{n_1}$ is decreased by $1$, the distance between $u$ and any vertex of $K_{n_2-1}$ is increased by $1$, and the distance between any other two vertices is unchanged.
Therefore,
$$
H(G_{n_1,n_2,n_3})-H(G_{n_1+1,n_2-1,n_3})  =  \left(\frac{n_1}{2}+n_2-1\right)-\left(n_1+\frac{n_2-1}{2}\right) =  \frac{n_2-n_1-1}{2} <  0.$$
The result follows.\hfill$\blacksquare$

%In the following, firstly, we determine the unique graph with the maximum Harary index among all the graphs with order $n$ and connectivity $r$.
\begin{theorem} \label{vercon}
For each $r=1,2,\ldots,n-2$, the graph $K(n-1,r)$ is the unique one with the maximum Harary index among all graphs of order $n$ and vertex connectivity $r$.
\end{theorem}

{\it Proof:}
Let $G$ be a graph that attains the maximum Harary index in $\mathcal{G}_n^r$.
Let $U$ be a vertex cut of $G$ containing $r$ vertices such that $G-U$ has components $G_1,G_2,\ldots,G_s$, where $s\geq 2$.
Firstly, we assert that $s=2$; otherwise adding all possible edges within the graph $G_1\cup G_2\cup\ldots \cup G_{s-1}$, we would get a graph belonging to $\mathcal{G}_n^r$ but with a larger Harary index.
Similarly, the induced subgraph $G[U]$, and the subgraphs $G_1,G_2$ are all complete, and each vertex of $U$ joins all vertices of $G_1$ and $G_2$.
Without loss of generality, we assume that $|V(G_1)|=:n_1 \geq |V(G_2)|=:n_2$.
By Lemma \ref{con}, $n_2=1$, and hence $G= K(n-1,r)$.\hfill$\blacksquare$

%Secondly, we determine the unique graph with the maximum Harary index among all the graphs with order $n$ and edge connectivity $r$.
\begin{theorem} \label{edgecon}
For each $r=1,2,\ldots,n-2$, the graph $K(n-1,r)$ is the unique one with the maximum Harary index among all graphs of order $n$ and edge connectivity $r$.
\end{theorem}

{\it Proof:}
Let $G$ be a graph that attains the maximum Harary index in $\overline{\mathcal{G}}_n^r$.
Assume the vertex connectivity of $G$ is $r_0$. Then $r_0 \leq r$.
So $$H(G) \leq H(K(n-1,r_0)) \le H(K(n-1,r)) \le H(G),$$
where the first inequality holds by Theorem \ref{vercon}, the second equality holds by Lemma \ref{addedge} as $K(n-1,r)$ is obtained from $K(n-1,r_0)$ 
by adding $r-r_0$ edges, and the last
inequality holds as $K(n-1,r) \in \overline{\mathcal{G}}_n^r$.
Hence, all inequalities above become equalities, which implies $r=r_0$, and $G=K(n-1,r)$ from the first equality by Theorem \ref{vercon}.
\hfill$\blacksquare$

\begin{coro}
Let $G$ be a graph of order $n$ with vertex or edge connectivity $r$, where $1 \le r \le n-2$.
Then
$$H(G)\leq \frac{(n-1)^2+r}{2},$$ with equality holds if and only if $G= K(n-1,r)$.
%For each $r=1,2,\ldots,n-2$ and $G \in \mathcal{G}_n^r(\overline{\mathcal{G}}_n^r,$respectively$)$, $
\end{coro}

{\it Proof:} By Theorems \ref{vercon} and \ref{edgecon}, we only need to calculate the Harary index of $K(n-1,r)$.
Since $K(n-1,r)$ is a graph with diameter $2$, the number of pairs of vertices with distance $1$ is $C_{n-1}^2+r$,
and the number of pairs of vertices with distance $2$ is $n-r-1$, we have
$$H(K(n-1,r))=C_{n-1}^2+r+\frac{n-r-1}{2}=\frac{(n-1)^2+r}{2}.$$
\hfill$\blacksquare$

Finally, we characterize the graphs with the second maximum Harary index among all the graphs of order $n$ with vertex connectivity $r$.
%Denote by $G_{n-r-1,1,r}-e_1,G_{n-r-1,1,r}-e_3,G_{n-r-1,1,r}-e_{13}$ the graphs obtained from $G_{n-r-1,1,r}$ by deleting an arbitrary edge within $K_{n-r-1}, K_r$ and between $K_{n-r-1}$ and $K_r$, respectively.

\begin{theorem}
Let $G$ be a graph with the second maximum Harary index among all graphs of order $n$ and vertex connectivity $r$, where $1 \le r \le n-2$.\\
{\em(1)} If $r=n-2$, then $G$ is obtained from $K_{n-2} \vee O_2$ by deleting an arbitrary edge in $K_{n-2}$.\\
{\em(2)} If $1\leq r \leq n-3$ and $r \neq n-4$, then $G$ is a graph obtained from $K(n-1,r)$ by deleting an arbitrary edge in the induced subgraph $K_{n-1}$.\\
{\em(2)} If $r=n-4$, then $G=G_{2,2,n-4}$ or $G$ is obtained from $K(n-1,r)$ by deleting an arbitrary edge in the induced subgraph $K_{n-1}$.
\end{theorem}

{\it Proof:} Let $U$ be a vertex cut of $G$ containing $r$ vertices such that $G-U$ has components $G_1,G_2,\ldots,G_s$,
where $|V(G_1)|\geq |V(G_2)|\geq \cdots \geq |V(G_s)|$,  and $s\geq 2$.
If $r=n-2$, then $s=2$, and $G$ is obtained from $K(n-1,n-2)=K_{n-2} \vee O_2$ by deleting an arbitrary edge in $K_{n-2}$.
So we assume $r \le n-3$.

We first claim that $s=2$, or $s=3$ and $|V(G_1)|= |V(G_2)|= |V(G_3)|=1$.
Otherwise, $s\geq 4$.  By connecting one vertex of $G_1$ with one of $G_2$, we will arrive at a new graph $G'$.
Obviously, $G' \in \mathcal{G}_n^r$ but $G'\neq K(n-1,r)$; and by Theorem \ref{vercon}, $H(K(n-1,r))>H(G')>H(G)$, a contradiction.
If $s=3$ and $|V(G_1)|\geq 2$. Then we attain a new graph $G''$ by connecting one vertex of $G_1$ with one of $G_2$.
Similarly we also have $H(K(n-1,r))>H(G'')>H(G)$, a contradiction.

If $s=3$ and $G_1, G_2, G_3$ are all single points, then $G = O_3\vee K_{n-3}$.
Now suppose $s=2$. By Lemma \ref{con}, 
if $|V(G_2)|=1$, then $G \in \{ G_{n-r-1,1,r}-e_1, G_{n-r-1,1,r}-e_3, G_{n-r-1,1,r}-e_{13}\}$;
if $|V(G_2)|>1$, then $r \leq n-4$ and $G= G_{n-r-2,2,r}$;
where $e_1,e_3$ are respectively the (arbitrary) edges in $K_{n-r-1}, K_{r}$ and $e_{13}$ is an (arbitrary) edges connecting $K_{n-r-1}$ and $K_{r}$,
by recalling $G_{n-r-1,1,r}=(K_{n-r-1} \cup K_1) \vee K_r$.
Observe that $O_3 \vee K_{n-3}= G_{n-r-1,1,r}-e_1$, where $r=n-3$.

By a little calculation,
$$
H(G_{n-r-1,1,r}-e_1)=H(G_{n-r-1,1,r}-e_3)=H(G_{n-r-1,1,r}-e_{13})=\frac{n^2-2n+r}{2}$$
and when $r\leq n-4$
$$H(G_{n-r-2,2,r})=\frac{n^2-3n+2r+4}{2}.$$
If $1\leq r<n-4$, $\frac{n^2-2n+r}{2}>\frac{n^2-3n+2r+4}{2}$; if $r=n-4$, $\frac{n^2-2n+r}{2}=\frac{n^2-3n+2r+4}{2}$.
So, if $1\leq r\leq n-3$ and $r\neq n-4$, then $G$ is one of $G_{n-r-1,1,r}-e_1,G_{n-r-1,1,r}-e_3$ or $G_{n-r-1,1,r}-e_{13}$, namely $G$ is obtained from $K(n-1,r)$ by deleting an arbitrary edge in the induced subgraph $K_{n-1}$.
If $r=n-4$, then $G=G_{2,2,n-4}$ or $G$ is obtained from $K(n-1,r)$ by deleting an arbitrary edge in the induced subgraph $K_{n-1}$.
\hfill$\blacksquare$

\end{document}